# ON $r$-HELIX HYPERSURFACES

**Evren ZIPLAR**

Department of Mathematics, Faculty of Science, University of Ankara, Tandoğan, Ankara,

**Ali ŞENOL**

Cankırı Karatekin University, Faculty of Science, Department of Mathematics, Cankırı,

**Yusuf YAYLI**

Department of Mathematics, Faculty of Science, University of Ankara, Tandoğan, Ankara,

**Abstract**

In this paper, we study strong $r$-helix hypersurfaces and the special curves on these surfaces. Moreover, we investigated the relations between strong $r$-helix hypersurfaces and the Gauss transformations of these surfaces in Euclidean $n$-space.

**Keywords:** Strong $r$ -helix hypersurface, Gauss transformation, line of curvature, geodesic curve, asymptotic curve, slant helix.

**Mathematics Subject Classification 2000:** 53A04, 53B25, 53C40, 53C50

## 1. Introduction

In differential geometry of surfaces, an helix hypersurface in $E^n$ is defined by the property that tangent planes make a constant angle with a fixed direction (helix direction) in [3]. Di Scala and Ruiz- Hernández have introduced the concept of these surfaces in [3]. Moreover, the concept of strong $r$-helix submanifold of $\mathrm{IR}^n$ was introduced in [2]. Let $M \subset \mathrm{IR}^n$ be a submanifold and let $H(M)$ be the set of helix directions of $M$. If $H(M)$ is $r$-dimensional linear subspace of $\mathrm{IR}^n$, then $M$ is called a strong $r$-helix [2].

Nowadays, M. Ghomi worked out the shadow problem given by H.Wente. And, He mentioned the shadow boundary in [8]. Ruiz- Hernández investigated that shadow bounderies are related to helix submanifolds whose tangent space makes constant angle with a fixed direction in [6].

Helix hypersurfaces has been worked in nonflat ambient spaces in [4,5]. Cermelli and Di Scala have also studied helix hypersurfaces in liquid cristals in [9].

A.I. Nistor has also introduced certain constant angle surfaces constructed on curves in $E^3$ in [1]. Özkaldi and Yaylı give some characterization for a curve lying on a surface for which the unit normal makes a constant angle with a fixed direction in [10].

One of the main purposes of this work is to observe the relations between strong $r$-helix hypersurfaces and special curves in Euclidean $n$-space $E^n$. Another purpose of this study is to give the relations between strong $r$-helix hypersurfaces and the Gauss transformations of these surfaces in Euclidean $n$-space $E^n$.

## 2. Preliminaries

**Definition 2.1** Let $\alpha : I \subset \mathrm{IR} \to E^n$ be an arbitrary curve in $E^n$. Recall that the curve $\alpha$ is said to be of unit speed ( or parametrized by the arc-length function $s$ ) if $\langle \alpha'(s), \alpha'(s) \rangle = 1$, where $\langle , \rangle$ is the standart scalar product in the Euclidean space $E^n$ given by

$$\langle X, Y \rangle = \sum_{i=1}^{n} x_i y_i,$$

for each $X = (x_1, x_2,..., x_n), Y = (y_1, y_2,..., y_n) \in E^n$.

Let $\{V_1(s), V_2(s),..., V_n(s)\}$ be the moving frame along $\alpha$, where the vectors $V_i$ are mutually orthogonal vectors satisfying $\langle V_i, V_i \rangle = 1$. The Frenet equations for $\alpha$ are given by

$$\begin{bmatrix} V_1' \\ V_2' \\ V_3' \\ \vdots \\ V_{n-1}' \\ V_n' \end{bmatrix} = \begin{bmatrix} 0 & k_1 & 0 & 0 & \cdots & 0 & 0 \\ -k_1 & 0 & k_2 & 0 & \cdots & 0 & 0 \\ 0 & -k_2 & 0 & k_3 & \cdots & 0 & 0 \\ \vdots & \vdots & \vdots & \vdots & & \vdots & \vdots \\ 0 & 0 & 0 & 0 & \cdots & 0 & k_{n-1} \\ 0 & 0 & 0 & 0 & \cdots & -k_{n-1} & 0 \end{bmatrix} \begin{bmatrix} V_1 \\ V_2 \\ V_3 \\ \vdots \\ V_{n-1} \\ V_n \end{bmatrix}.$$

Recall that the functions $k_i(s)$ are called the $i$-th curvatures of $\alpha$ [7].

**Definition 2.2** Let $\alpha : I \subset \mathrm{IR} \to E^n$ be a unit speed curve with nonzero curvatures $k_i$ ($i = 1, 2, ..., n$) in $E^n$ and let $\{V_1, V_2, ... V_n\}$ denote the Frenet frame of the curve $\alpha$. We call $\alpha$ a $V_n$-slant helix, if the $n$-th unit vector field $V_n$ makes a constant angle $\varphi$ with a fixed direction X, that is,

$$\langle V_n, \mathrm{X} \rangle = \cos(\varphi),\ \varphi \neq \frac{\pi}{2},\ \varphi = \text{constant}$$

along the curve, where X is unit vector field in $E^n$ [7].

**Definition 2.3** Given a hypersurface $M \subset \mathrm{IR}^n$ and an unitary vector $d \neq 0$ in $\mathrm{IR}^n$, we say that $M$ is a helix with respect to the fixed direction $d$ if for each $q \in M$ the angle between $d$ and $T_q M$ is constant. Note that the above definition is equivalent to the fact that $\langle d, \xi \rangle$ is constant function along $M$, where $\xi$ is a unit normal vector field on $M$ [3].

**Definition 2.4** A submanifold $M \subset \mathrm{IR}^n$ is a $r$-helix if there is exist a linear subspace $H \subset \mathrm{IR}^n$ of dimension $r = \dim(H)$ such that $M$ is a helix with respect to any direction $d \in H$. The subspace $H$ is called the subspace of helix directions [3].

**Definition 2.5** Let $M \subset \mathrm{IR}^n$ be a submanifold of a euclidean space. A vector $d \in \mathrm{IR}^n$ is called a helix direction of $M$ if the angle between $d$ and any tangent space $T_p M$ is constant. Let $H(M)$ be the set of helix directions of $M$. We say that $M$ is a strong $r$-helix if $H(M)$ is $r$-dimensional linear subspace of $\mathrm{IR}^n$ [2].

**Theorem 2.1** Let $\alpha : I \subset \mathrm{IR} \to E^n$ be a unit speed curve ( parametrized by arclength function $s$) in $E^n$ and let $\beta : I \to S^{n-1} \subset E^n$ be the tangent indicatrix of the curve $\alpha$, where $S^{n-1}$ is the unit hypersphere in $E^n$. Then the curve $\alpha$ is a slant helix with direction $L$ in $E^n$ if and only if the curve $\beta$ is a general helix (spherical helix) with direction $L$ on $S^{n-1} \subset E^n$. In other words, $\alpha$ and $\beta$ have the same direction $L$ [11].

**Theorem 2.2** Let $M$ be a $r-$helix hypersurface in $E^n$ and let $H$ be the subspace of helix directions of $M$. If $\alpha : I \subset \mathbb{R} \to M$ be a unit speed geodesic curve on $M$, then the curve $\alpha$ is a $V_2$-slant helix with respect to any direction $d \in H$ in $E^n$ [12].

**Proof:** Let $\xi$ be a normal vector field on $M$. Since $M$ is $r-$ helix hypersurface, $<d, \xi> = $ constant, where $d \in H$ is any direction. That is, the angle between $d$ and $\xi$ is constant on every point of the surface $M$. And, $\alpha''(s) = \lambda \xi \,|_{\alpha(s)}$ along the curve $\alpha$ since $\alpha$ is a geodesic curve on $M$. Moreover, by using the Frenet equation $\alpha''(s) = V_1' = k_1 V_2$, we obtain $\lambda \xi \,|_{\alpha(s)} = k_1 V_2$, where $k_1$ is first curvature of $\alpha$. Thus, from the last equation, by taking norms on both sides, we obtain $\xi = V_2$ or $\xi = -V_2$. So, $\langle d, V_2 \rangle$ is constant along the curve $\alpha$ since $\langle d, \xi \rangle = $ constant. In other words, the angle between $d$ and $V_2$ is constant along the curve $\alpha$. Consequently, the curve $\alpha$ is a $V_2$-slant helix with respect to any direction $d \in H$ in $E^n$.

## 3. MAIN THEOREMS

**Theorem 3.1** Let $\alpha : I \subset \mathbb{R} \to M \subset E^n$ be a geodesic curve on $M$ with unit speed (or parametrized by arclength function $s$) and let $H \subset \mathbb{R}^n$ be the subspace of helix directions of $M$, where $M$ is a $r$-helix hypersurface. Then the tangent indicatrix $\alpha'(s)$ of the curve $\alpha(s)$ is a spherical helix (general helix) with respect to any helix direction in $H$ on unit the hypersphere $S^{n-1} \subset E^n$.

**Proof:** We assume that $\alpha$ is a geodesic curve on $M$. Then from Theorem 2.2 the curve $\alpha$ is a $V_2$-slant helix with respect to any helix direction $d \in H$ in $E^n$. On the other hand, from Theorem 2.1, the tangent indicatrix $\alpha'$ of the curve $\alpha$ with the direction $d$ is a spherical helix with respect to the same direction $d$ in $H$ on the unit hypersphere $S^{n-1} \subset E^n$. Consequently, $\alpha'$ is a spherical helix (general helix) with respect to any helix direction in $H$ on the unit hypersphere $S^{n-1} \subset E^n$. This completes the proof.

**Theorem 3.2** Let $M$ be a hypersurface in Euclidean $n$-space $E^n$ and let $N$ be the unit normal vector field of $M$. If $M$ is a strong $r$-helix in $E^n$, where $H(M) \subset \mathrm{IR}^n$ is the space of helix directions of $M$, then the position vectors on all points of the surface
$$\eta(M) = \{X \in S^{n-1} \mid X = N(P), P \in M\}$$
make a constant angle with the space $H(M) \subset \mathrm{IR}^n$. Here,
$$\eta : M \to S^{n-1} \subset E^n$$
$$P \to \eta(P) = \vec{N}(P)$$
is Gauss transformation of $M$, where $S^{n-1}$ is the unit hypersphere in $E^n$.

**Proof:** Let $d \in H(M)$ be any helix direction of $M$. Since $M$ is a strong $r$-helix in $E^n$,

$$< N(P), d > = \text{constant}$$

for every $P \in M$. So, we obtain

$$< X, d > = \text{constant},$$

for every $X = \vec{N}(P) \in \eta(M)$, where $d \in H(M)$ any helix direction of $M$. It follows that
$$< X, d > = \text{constant}$$
for all $d \in H(M)$, where $X$ is the position vector of $\mu(M)$. Hence, the position vectors on all points of the surface $\mu(M)$ make a constant angle with the space $H(M) \subset \mathrm{IR}^n$. This completes the proof.

**Theorem 3.3** Let $M$ be a strong $r$-helix hypersurface in Euclidean $n$-space $E^n$ and let $H(M) \subset \mathrm{IR}^n$ be the space of helix directions of $M$. If $\alpha : I \subset \mathrm{IR} \to M$ $(\alpha(t) \in I, t \in I)$ is a curve on the surface $M$ and $\exists d_j \notin T_p M$ for every $p = \alpha(t) \in M$, where $d_j \in H(M)$, then the curve

$$\eta(\alpha(t)) = \{X \in S^{n-1} \mid X = N\mid_{\alpha(t)}, \alpha(t) \in M\}$$

is not a geodesic curve on the unit hypersphere $S^{n-1} \subset E^n$, where $N\mid_{\alpha(t)}$ the unit normal vector field of $M$ along the curve $\alpha$ and

$$\eta : M \to S^{n-1} \subset \mathrm{E}^n$$
$$\alpha(t) \to \eta(\alpha(t)) = \vec{N}\,|_{\alpha(t)}$$

Gauss transformation of $M$ on the curve $\alpha$.

**Proof:** Since $M$ is a strong $r$-helix hypersurface, for all $d_j \in H(M)$,
$$< N\,|_{\alpha(t)}, d_j > = \cos(\theta_j) = \text{constant}$$

along the curve $\alpha$. We assume that $\beta = \eta(\alpha(t))$, where $\eta$ Gauss transformation of $M$. So, we obtain

$$\beta = N\,|_{\alpha(t)}$$

along the curve $\alpha$ since $\eta(\alpha(t)) = N\,|_{\alpha(t)}$ along the curve. Thus, we have

$$< \beta, d_j > = \cos(\theta_j) = \text{constant}$$

for all $d_j \in H(M)$.

If we take the derivative in each part of the equality $< \beta, d_j > = \cos(\theta_j)$ twice, we obtain, for all $d_j \in H(M)$,

$$< \beta'', d_j > = 0. \qquad (1)$$

Now, we suppose that the curve $\beta$ is a geodesic curve on $S^{n-1}$. Then $\beta'' = kN\,|_{\alpha(t)}$, where $N\,|_{\alpha(t)}$ is also the normal vector field of $S^{n-1}$. Hence, we get, for all $d_j \in H(M)$,

$$< N\,|_{\alpha(t)}, d_j > = 0$$

by using the equation (1). It follows that all $d_j \in T_{\alpha(t)}M$ for every point $\alpha(t) \in M$. But, according to the hypothesis in this Theorem, $\exists d_j \notin T_p M$ for every $p = \alpha(t) \in M$. So, it is a contradiction. As a result, $\beta = \eta(\alpha(t))$ is not a geodesic curve on the unit hypersphere $S^{n-1} \subset \mathrm{E}^n$.

**Theorem 3.4** Let $M$ be strong $r-$ helix hypersurface in $E^n$ and $H(M) \subset \mathrm{IR}^n$ be the space of helix directions of $M$. If $\alpha : \mathrm{I} \subset \mathrm{IR} \to M$ ( $\alpha(t) \in M$, $t \in \mathrm{I}$ ) is a curve on the surface $M$, then

$d \in Sp\{\beta'\}^{\perp}$ along the curve $\beta$ for all $d \in H(M)$. Here,

$$\beta = \eta(\alpha(t)) = \{X \in S^{n-1} \mid X = N\mid_{\alpha(t)}, \alpha(t) \in M\}$$

and

$$\eta: M \to S^{n-1} \subset E^n$$
$$\alpha(t) \to \eta(\alpha(t)) = \vec{N}\mid_{\alpha(t)}$$

is Gauss transformation of $M$ on the curve $\alpha$, where $S^{n-1}$ is the unit hypersphere in $E^n$ and $N\mid_{\alpha(t)}$ the unit normal vector field of $M$ along the curve $\alpha$.

**Proof:** Since $M$ is a strong $r$−helix hypersurface in $E^n$,

$$< N\mid_{\alpha(t)}, d > = \text{constant}$$

along the curve $\alpha$ for all $d \in H(M)$. If we take the derivative in each part of the equality $< N\mid_{\alpha(t)}, d > = \text{constant}$, we have

$$< N'\mid_{\alpha(t)}, d > = 0.$$

On the other hand, $\beta' = N'\mid_{\alpha(t)}$ since $\beta = N\mid_{\alpha(t)}$. So, we obtain:

$$< \beta', d > = 0$$

throughout the curve $\beta$. As a result, $d \in Sp\{\beta'\}^{\perp}$ along the curve $\beta$ for all $d \in H(M)$.

**Corollary 3.1** Let $M$ be helix hypersurface with the fixed direction $d$ in $E^n$. If $\alpha: I \subset \mathbb{R} \to M$ ($\alpha(t) \in M, t \in I$) is a curve on the surface $M$, then $d \in Sp\{\beta'\}^{\perp}$ along the curve $\beta$. Here,

$$\beta = \eta(\alpha(t)) = \{X \in S^{n-1} \mid X = N\mid_{\alpha(t)}, \alpha(t) \in M\}$$

and

$$\eta : M \to S^{n-1} \subset E^n$$
$$\alpha(t) \to \eta(\alpha(t)) = \vec{N}|_{\alpha(t)}$$

is Gauss transformation of $M$ on the curve $\alpha$, where $S^{n-1}$ is the unit hypersphere in $E^n$ and $N|_{\alpha(t)}$ the unit normal vector field of $M$ along the curve $\alpha$.

**Theorem 3.5** Let $M$ be strong $r-$helix hypersurface in $E^n$ and let $H(M)$ be the space of helix directions of $M$. If $\alpha : I \subset \mathrm{IR} \to M$ ($\alpha(s) \in M$, $s \in I$) is a curve with unit speed on the surface $M$ ($s$ arc-length parameter) and if $\alpha$ is a geodesic curve on $M$, then $d \in Sp\{V_2'\}^\perp$ along the curve $\alpha$, where $V_2$ is the element of the moving frame $\{V_1(s), V_2(s), ..., V_n(s)\}$ of the curve $\alpha$ and $d \in H(M)$ is any direction.

**Proof:** Since $M$ is strong $r-$ helix hypersurface, for any direction $d \in H(M)$,

$$< N|_{\alpha(s)}, d > = \text{constant}$$

along the curve $\alpha$, where $N$ is the normal vector field of $M$. If we are taking the derivative in each part of the equality with respect to $s$, we obtain :

$$\langle (N|_{\alpha(s)})', d \rangle = 0 \qquad (2)$$

along the curve $\alpha$. And, $\alpha''(s) = \lambda N|_{\alpha(s)}$ along the curve $\alpha$ since $\alpha$ is a geodesic curve on $M$. Moreover, by using the Frenet equation $\alpha''(s) = V_1' = k_1 V_2$, we obtain $\lambda N|_{\alpha(s)} = k_1 V_2$, where $k_1$ is first curvature of $\alpha$. Thus, from the last equation, by taking norms on both sides, we obtain $N|_{\alpha(s)} = V_2$ or $N|_{\alpha(s)} = -V_2$. So, $(N|_{\alpha(s)})' = V_2'$. It follows that $< V_2', d > = 0$ by using the equality (2). Consequently, $d \in Sp\{V_2'\}^\perp$ along the curve $\alpha$.

**Corollary 3.2** Let $M$ be a helix hypersurface with the direction $d$ in $E^n$ and let $\alpha : I \subset \mathrm{IR} \to M$ ($\alpha(s) \in M$, $s \in I$) be a curve with unit speed on the surface $M$ (s arc-length parameter). If $\alpha$ is a geodesic curve on $M$, then $d \in Sp\{V_2'\}^{\perp}$ along the curve $\alpha$, where $V_2$ is the element of the moving frame $\{V_1(s), V_2(s), ..., V_n(s)\}$ of the curve $\alpha$.

**Theorem 3.6** Let $M$ be a strong $r$-helix hypersurface in $E^n$ and let $H(M) \subset \mathrm{IR}^n$ be the space of the helix directions of $M$. Let $\alpha: I \subset \mathrm{IR} \to M$ ($\alpha(s) \in M$, $s \in I$) be a curve with unit speed on the surface $M$ (s arc-length parameter). Let us assume that for each $q \in M$ the angle $\theta_j \neq \{0, \pi/2\}$ between each $d_j \in H(M)$ and $N$. If for $\exists d_j \in H(M)$ $d_j \in Sp\{V_1, N\}$ along the curve $\alpha$, then the curve $\alpha$ is a asymptotic curve on the surface $M$, where $V_1$ is the element of the Frenet frame $\{V_1, V_2, ..., V_n\}$ and $N$ is the unit normal vector field of the surface $M$.

**Proof:** Since $M$ is a strong $r$-helix hypersurface and $d_j \in Sp\{V_1, N\}$ for $\exists d_j \in H(M)$ along the curve $\alpha$, we can write:
$$d_j = \cos(\theta_j) N + \sin(\theta_j) V_1. \tag{3}$$
Taking the derivative in each part of the equation (3), we obtain:
$$0 = \cos(\theta_j) N' + \sin(\theta_j) V_1' \tag{4}$$
And, doing the scalar product with $V_1$ in each part of the equation (4), we have:
$$0 = \cos(\theta_j) <N', V_1> + \sin(\theta_j) <V_1', V_1> \tag{5}$$
On the other hand, $<V_1', V_1> = 0$ since $V_1$ is a unit vector. So, we get from the equation (5):
$$\cos(\theta_j) <N', V_1> = 0.$$
According to hypothesis in this Theorem, since $\cos(\theta_j) \neq 0$, it follows that
$$<N', V_1> = 0.$$
Finally, the curve $\alpha$ is a asymptotic curve. This completes the proof.

**Corollary 3.3** Let $M$ be a strong $r$-helix hypersurface in $E^n$ and let $H(M) \subset \mathrm{IR}^n$ be the space of the helix directions of $M$. Let $\alpha: I \subset \mathrm{IR} \to M$ ($\alpha(s) \in M$, $s \in I$) be a geodesic curve with unit speed on the surface $M$ (s arc-length parameter). Let us assume that for each $q \in M$ the angle $\theta_j \neq \{0, \pi/2\}$ between each $d_j \in H(M)$ and $N$. If for $\exists d_j \in H(M)$ $d_j \in Sp\{V_1, N\}$ along the curve $\alpha$, then the curve $\alpha$ is a line on the surface $M$, where $V_1$ is the element of the Frenet frame $\{V_1, V_2, ..., V_n\}$ and $N$ is the unit normal vector field of the surface $M$.

**Proof:** From Theorem 3.6, the curve $\alpha$ is a asymptotic curve. Moreover, according to the hypothesis in this Theorem, the curve $\alpha$ is a geodesic curve. As we know, if a curve is both asymptotic and geodesic, then the curve is a line only.

**Corollary 3.4** Let $M$ be a strong $r$-helix hypersurface in $\mathrm{E}^n$ and let $H(M) \subset \mathrm{IR}^n$ be the space of the helix directions of $M$. Let $\alpha : \mathrm{I} \subset \mathrm{IR} \to M$ ($\alpha(s) \in M$, $s \in \mathrm{I}$) be a line of curvature with unit speed on the surface $M$ (s arc-length parameter). Let us assume that for each $q \in M$ the angle $\theta_j \neq \{0, \pi/2\}$ between each $d_j \in H(M)$ and $N$. If for $\exists d_j \in H(M)$ $d_j \in Sp\{V_1, N\}$ along the curve $\alpha$, then the curve $\alpha$ is a line on the surface $M$, where $V_1$ is the element of the Frenet frame $\{V_1, V_2, ..., V_n\}$ and $N$ is the unit normal vector field of the surface $M$.

**Proof:** From Theorem 3.6, the curve $\alpha$ is a asymptotic curve. Moreover, according to the hypothesis in this Theorem, the curve $\alpha$ is a line of curvature. As we know, if a curve is both asymptotic and line of curvature, then the curve is a line only.